\documentclass[twoside,openany]{book}
\usepackage{leftidx}
\usepackage{graphicx}
\usepackage{amssymb}
\usepackage{mathrsfs}
\usepackage{multirow}
\usepackage{amsmath,hyperref}
\usepackage{amsfonts,amsthm}
\usepackage{fancyhdr}
\usepackage{verbatim}
\usepackage{slashbox}
\usepackage{rotating}
\usepackage{booktabs}
\usepackage[stable]{footmisc}
\usepackage{indentfirst,cite}
\usepackage[all]{xy}
\usepackage{amscd}
\usepackage{rotating}
\usepackage[sectionbib]{chapterbib}
\usepackage{makeidx}\makeindex
\evensidemargin=0pt \DeclareSymbolFont{lettersA}{U}{txmia}{m}{it}
\DeclareMathSymbol{\piup}{\mathord}{lettersA}{25}
\DeclareMathSymbol{\muup}{\mathord}{lettersA}{22}

\newsymbol\wjzhml 203F%空集

\newtheorem{theorem}{Theorem}

\renewcommand{\leq}{\leqslant}
\renewcommand{\geq}{\geqslant}

\theoremstyle{definition}

\newtheorem{example}[theorem]{Example}

\newtheorem*{ques0}{Question}

\numberwithin{equation}{section} \numberwithin{figure}{section}
\renewcommand{\thefigure}{\thesection.\arabic{figure}}
\begin{document}
\def\theequation{\thechapter.\arabic{equation}}\setcounter{equation}{0}

\setcounter{chapter}{1}

\baselineskip 16pt\renewcommand{\baselinestretch}{1.25}

\def\<{\langle}
\def\>{\rangle}
\def\({\left(}
\def\){\right)}
\def\[{\left[}
\def\]{\right]}

\newpage
\vspace*{11mm}\thispagestyle{empty}
\begin{center}{\LARGE\bf Recent Progress on Austere Submanifolds%\\[3mm]
}\end{center}
\thispagestyle{empty}\setcounter{page}{1}\vspace{5mm}

\markboth{\small {\sc Jianquan Ge \& Yi Zhou}}{\small Recent Progress on Austere Submanifolds}

\footnotetext{J. Q. Ge is partially supported by NSFC (No. 12171037) and the Fundamental Research Funds for the Central Universities.}
\footnotetext{Y. Zhou is partially supported by  NSFC (No. 12171037, 12271040), China Postdoctoral Science Foundation (No. BX20230018)
and National Key R$\&$D Program of China 2020YFA0712800.}

\noindent{\large {\sc Jianquan Ge}}

\noindent {\small School of Mathematical Sciences, Laboratory of Mathematics and Complex Systems, Beijing Normal University, Beijing, China}

\noindent {\small E-mail: jqge@bnu.edu.cn}\vspace{3mm}

\noindent{\large {\sc Yi Zhou}}

\noindent {\small Beijing International Center for Mathematical Research, Peking University, Beijing, China}

\noindent {\small E-mail: yizhou@bicmr.pku.edu.cn}\vspace{3mm}

\begin{quotation}
\noindent {\bf Abstract}\quad
{\small This is a brief survey of recent results related to austere submanifolds,
mainly based on the papers \cite{GTZ24, GZ23}.}

\end{quotation}

\section*{1. Backgrounds}\setcounter{equation}{0}\setcounter{chapter}{1}

A submanifold $M^n$ of a Riemannian manifold is said to be austere if
all symmetric polynomials of odd degree in the principal curvatures
with respect to any normal vector vanish, or equivalently,
the nonzero principal curvatures in every normal direction occur in oppositely signed pairs.
This concept was introduced by Harvey-Lawson \cite{HL} for constructing special Lagrangian submanifolds in $\mathbb{C}^N$.
It is not hard to see that austerity implies minimality,
and these two concepts are equivalent if and only if $n=1$ or $2$.

In the celebrated paper \cite{Bryant}, Bryant developed some algebra to describe the possible second fundamental forms of austere submanifolds.
Let $\operatorname{II}$ denote the second fundamental form,
and for any $p\in M$, let $|\operatorname{II}_p|\subset S^2(T_p^*M)$ be the subspace spanned by
$\{\operatorname{II}_{\xi}: \xi\in T_p^\bot M\}$,
where $\operatorname{II}_{\xi}(X, Y):=\<\operatorname{II}(X, Y),\xi\>$ is the second fundamental form
along the normal direction $\xi$.
By fixing an orthonormal basis, we can regard the vector space $|\operatorname{II}_p|$
as a subspace of the vector space $\mathcal{S}_n$ of $n\times n$ real symmetric matrices.
For convenience, we say that an $n\times n$ real symmetric matrix $A$ is austere if
its nonzero eigenvalues occur in oppositely signed pairs,
or equivalently, $$\operatorname{tr}A^{2k+1}=0\ \mbox{for each}\ 0\leq k\leq \[\frac{n-1}2\].$$
Thus the problem is translated to classify the subspaces consisting of austere matrices, called austere subspaces.
Note that the sum of two austere matrices is not necessary an austere matrix.
In order to construct examples of austere subspaces,
there is an insightful observation that for any orthogonal matrix $P\in\operatorname{O}(n)$,
the subspace $\mathcal{Q}_P:=\{A\in\mathcal{S}_n: PA+AP=0\}$ is austere.
We then present two classes of concrete examples.
\begin{example}
If $n$ is even and $P^2=-I_n$, then up to a conjugate action of $\operatorname{O}(n)$, we can take
$P=\begin{pmatrix}
0 & I_{\frac{n}2}  \\
-I_{\frac{n}2}  & 0 \\
\end{pmatrix}$,
and thus the austere subspace $\mathcal{Q}_P$ is
$$\mathcal{Q}(n):=\left\{\begin{pmatrix}
a_1 & a_2  \\
a_2 & -a_1 \\
\end{pmatrix}: a_1, a_2\in\mathcal{S}_{\frac{n}2}\right\}.$$
It is easy to see that
almost complex submanifolds of K\"{a}hler manifolds realize this class of austere subspaces.
\end{example}

\begin{example}
If $P^2=I_n$ and $P\neq\pm I_n$, then there exist positive integers $p,q$ such that $P$ conjugates to
$\begin{pmatrix}
I_p & 0  \\
0  & -I_q \\
\end{pmatrix}$,
and the austere subspace $\mathcal{Q}_P$ conjugates to
$$\mathcal{Q}(p,q):=\left\{\begin{pmatrix}
0 & a  \\
a^t & 0 \\
\end{pmatrix}: a\in M(p, q)\right\},$$
where $M(p, q)$ denotes the vector space of $p\times q$ real matrices.
\end{example}

Through the above observation and some complicated discussions,
Bryant constructed several maximal austere subspaces (see Table \ref{max_austere})
and completed the classification for $n\leq4$ (see Table \ref{classaust}), up to the conjugate action of $\operatorname{O}(n)$ on $\mathcal{S}_{n}$ given by $P\cdot A=P^tAP$.

\begin{table}[ht!]\small
\caption{Bryant's constructions of maximal austere subspaces}\label{max_austere}
\centering
\begin{tabular}{|c|c|c|}
\hline
Order & Maximal austere subspaces & Dimension\\
\hline
\ & $\left\{\begin{pmatrix}
m_1 & m_2  \\
m_2 & -m_1 \\
\end{pmatrix}: m_1, m_2\in\mathcal{S}_{p+1}\right\}$ & $(p+1)(p+2)$\\
\cline{2-3}
$n=2p+2$  & $\left\{\begin{pmatrix}
0 & m  \\
m^t & 0 \\
\end{pmatrix}: m\in M(k, n-k, \mathbb{R})\right\},\ 1\leq k\leq p$ & $k(n-k)$\\
\cline{2-3}
\  & $\left\{\begin{pmatrix}
\lambda I_{p+1} & m  \\
m^t & -\lambda I_{p+1} \\
\end{pmatrix}: m\in M(p+1, \mathbb{R}),\ \lambda\in\mathbb{R}\right\}$ & $(p+1)^2+1$\\
\hline
\ & $\left\{\begin{pmatrix}
m_1 & m_2 & 0 \\
m_2 & -m_1 & 0 \\
0 & 0 & 0 \\
\end{pmatrix}: m_1, m_2\in\mathcal{S}_{p}\right\}$ & $p(p+1)$\\
\cline{2-3}
$n=2p+1$ & $\left\{\begin{pmatrix}
0 & m  \\
m^t & 0 \\
\end{pmatrix}: m\in M(k, n-k, \mathbb{R})\right\},\ 1\leq k\leq p$ & $k(n-k)$\\
\cline{2-3}
\ & $\left\{\begin{pmatrix}
\lambda I_{p} & m & 0 \\
m^t & -\lambda I_{p} & 0 \\
0 & 0 & 0 \\
\end{pmatrix}: m\in M(p, \mathbb{R}),\ \lambda\in\mathbb{R}\right\}$ & $p^2+1$\\
\hline
\end{tabular}
\end{table}

\begin{table}[ht!]\small
\caption{Bryant's classifications for $n=2,3,4$}\label{classaust}
\centering
\begin{tabular}{|c|c|c|}
\hline
Order & Maximal austere subspaces & Dimension\\
\hline
$n=2$ & $\left\{\begin{pmatrix}
a & b  \\
b & -a \\
\end{pmatrix}: a, b\in\mathbb{R}\right\}$ & $2$\\
\hline
$n=3$ & $\left\{\begin{pmatrix}
a & b & 0 \\
b & -a & 0 \\
0 & 0 & 0 \\
\end{pmatrix}: a, b\in\mathbb{R}\right\}$ & $2$\\
\cline{2-3}
\ & $\left\{\begin{pmatrix}
\,0 & 0 & a\, \\
\,0 & 0 & b\, \\
\,a & b & 0\, \\
\end{pmatrix}: a, b\in\mathbb{R}\right\}$ & $2$\\
\hline
$n=4$ & $\left\{\begin{pmatrix}
m_1 & m_2  \\
m_2 & -m_1 \\
\end{pmatrix}: m_1, m_2\in\mathcal{S}_{2}\right\}$ & $6$\\
\cline{2-3}
\  & $\left\{\begin{pmatrix}
\lambda I_{2} & m  \\
m^t & -\lambda I_{2} \\
\end{pmatrix}: m\in M(2, \mathbb{R}),\ \lambda\in\mathbb{R}\right\}$ & $5$\\
\cline{2-3}
\  & $\left\{\begin{pmatrix}
0 & x_1 & x_2 & x_3  \\
x_1 & 0 & \lambda_3x_3 & \lambda_2x_2 \\
x_2 & \lambda_3x_3 & 0 & \lambda_1x_1 \\
x_3 & \lambda_2x_2 & \lambda_1x_1 & 0 \\
\end{pmatrix}: x_1, x_2, x_3\in\mathbb{R}\right\}$, where the & $3$ \\
\  & constants $\lambda_1\geq\lambda_2\geq0\geq\lambda_3$ satisfy $\lambda_1\lambda_2\lambda_3+\lambda_1+\lambda_2+\lambda_3=0.$ & \ \\
\hline
\end{tabular}
\end{table}

Based on this classification of maximal austere spaces, Bryant locally classified $3$-dimensional austere submanifolds in Euclidean spaces in the same paper.
Then Ionel-Ivey \cite{II10, II12} continued to study $4$-dimensional austere submanifolds in Euclidean spaces and obtained some classification results under additional assumptions.
In addition, Dajczer-Florit \cite{DF01} and Choi-Lu \cite{CL08} locally classified
$3$-dimensional austere submanifolds in unit spheres and hyperbolic spaces, respectively.
More constructions and classification results for austere submanifolds can be found in \cite{DV23, GZ23, IST09, KM22}.

A rich source of austere submanifolds is the focal submanifolds of isoparametric hypersurfaces in unit spheres which are hypersurfaces of constant principal curvatures.
In \cite{Cartan1, Cartan2}, Cartan proved that an isoparametric family of hypersurfaces
in unit spheres with 3 distinct constant principal curvatures of equal multiplicity $m$ is given by
the regular level sets of the Cartan-M\"{u}nzner polynomials
\begin{equation}\label{Cartan-poly}
\begin{aligned}
F_C(x, y, X, Y, Z)=&x^3-3xy^2+\frac32x(X\overline{X}+Y\overline{Y}-2Z\overline{Z})\\
&+\frac{3\sqrt{3}}2y(X\overline{X}-Y\overline{Y})+\frac{3\sqrt{3}}2(XYZ+\overline{Z}\overline{Y}\overline{X}),
\end{aligned}
\end{equation}
where $x, y\in\mathbb{R}$ and $X, Y, Z\in\mathbb{R}, \mathbb{C}, \mathbb{H}, \mathbb{O}$ for
$m=1, 2, 4, 8$, respectively.
In fact, all these isoparametric hypersurfaces are extrinsically homogeneous,
that is, they are principal orbits in the unit spheres under the group actions of
$\operatorname{SO}(3)$, $\operatorname{SU}(3)$, $\operatorname{Sp}(3)$ or the exceptional Lie group $F_4$.
Moreover, Cecil-Ryan \cite{CR85, CR15} pointed out that the former three actions can be presented as conjugate actions on the vector space of $3\times3$ traceless real, complex Hermitian, or quaternionic Hermitian matrices.
For example, let $E^5$ be the 5-dimensional vector space of $3\times3$ real symmetric matrices with trace zero, and let $S^4$ be the unit sphere in $E^5$ with respect to the Frobenius norm.
Then the Cartan-M\"{u}nzner polynomial of (\ref{Cartan-poly}) is nothing but the determinant function (cf. \cite{GT22}), and the principal orbits of the isometric $\operatorname{SO}(3)$-action
$$\operatorname{SO}(3)\times S^4\rightarrow S^4,\ (P, A)\mapsto P^tAP,$$
are exactly the isoparametric hypersurfaces with 3 distinct principal curvatures of multiplicity $m=1$.
In particular, the unique minimal hypersurface in this isoparametric family
has constant principal curvatures $-\sqrt{3}, 0, \sqrt{3}$ with respect to any unit normal vector,
and thus it is an austere submanifold.

Let $\mathcal{A}_n$ denote the set of $n\times n$ austere matrices with Frobenius norm $1$.
We observed that the preceding minimal isoparametric hypersurface is the orbit of the matrix
$$\frac1{\sqrt2}\begin{pmatrix}
1 & 0 & 0 \\
0 & -1 & 0 \\
0 & 0 & 0 \\
\end{pmatrix},$$ and is exactly the set $\mathcal{A}_3$.
Hence the intersection of $S^4$ and every austere subspace of $\mathcal{S}_3$
should be a totally geodesic sphere contained in $\mathcal{A}_3$.
Motivated by these interesting observations, we want to generalize the submanifold structure of $\mathcal{A}_3$ to (a subset of) $\mathcal{A}_n$.

\section*{2. The Construction of New Examples}
In this section, we briefly introduce our construction in \cite{GZ23}
which can handle the complex and quaternionic versions simultaneously.
Let $\mathbb{F}$ be one of $\mathbb{R}, \mathbb{C}$ or $\mathbb{H}$,
and let $m:=\operatorname{dim}_{\mathbb{R}}(\mathbb{F})=1, 2$ or $4$, respectively.
Then the set of traceless Hermitian matrices over $\mathbb{F}$, denoted by
$$E(n, \mathbb{F}):=\{A\in M(n, \mathbb{F}): A^*=A, \operatorname{tr}A=0\},$$
is a $N(n, \mathbb{F})$-dimensional real subspace of $M(n, \mathbb{F})$,
where $A^*$ denotes the conjugate transpose of $A$ and
$$N(n, \mathbb{F}):=\frac12 n(n-1)m+n-1
=\left\{\begin{aligned}
\frac12 n(n+1)-1\ \ \ &\mbox{if}\ \mathbb{F}=\mathbb{R};\\
n^2-1\ \ \ \ \ \ \ \ &\mbox{if}\ \mathbb{F}=\mathbb{C};\\
2n^2-n-1\ \ \ \ \, &\mbox{if}\ \mathbb{F}=\mathbb{H}.
\end{aligned}
\right.$$
Let $S(n, \mathbb{F})$ be the unit sphere in $E(n, \mathbb{F})$ with respect to
the Frobenius norm.
%usual Euclidean inner product given by $\<A, B\>:=\mathfrak{R}(\operatorname{tr}AB^*)$,
%where $\mathfrak{R}$ denotes the real part.
For each $k\geq2$, let $$f_k: S(n, \mathbb{F})\rightarrow\mathbb{R},\ A\mapsto\operatorname{tr}A^k,$$
and let $$\Phi_{n, \mathbb{F}}: S(n, \mathbb{F})\rightarrow\mathbb{R}^p,\
A\mapsto(f_3(A), f_5(A), \cdots, f_{2p+1}(A)),$$
where $p=\[\frac{n-1}2\]$.
Then we have $$\Phi_{n, \mathbb{R}}^{-1}(0)=\left\{A\in S(n, \mathbb{R}):
f_{2k+1}(A)=0\ \mbox{for each}\ 0\leq k\leq \[\frac{n-1}2\]\right\}
=\mathcal{A}_n.$$
Let $\mathcal{B}_{n, \mathbb{F}}$ and $\mathcal{C}_{n, \mathbb{F}}$ denote
the subsets of regular points and critical points of $\Phi_{n, \mathbb{F}}$
in $\Phi_{n, \mathbb{F}}^{-1}(0)$, respectively.
We can show that $\mathcal{C}_{n, \mathbb{F}}$ is closed in $S(n, \mathbb{F})$,
and thus by the regular level set theorem, we know
$\mathcal{B}_{n, \mathbb{F}}$ is a $(N(n, \mathbb{F})-p-1)$-dimensional properly embedded submanifold of $S(n, \mathbb{F})\setminus\mathcal{C}_{n, \mathbb{F}}$.

Let $U(n, \mathbb{F})$ denote the Lie group
$\operatorname{SO}(n)$, $\operatorname{U}(n)$, $\operatorname{Sp}(n)$
for $\mathbb{F}=\mathbb{R}$, $\mathbb{C}$, $\mathbb{H}$, respectively.
To study basic properties of the submanifold $\mathcal{B}_{n, \mathbb{F}}$,
we need to consider the isometric action of $U(n, \mathbb{F})$ on $E(n, \mathbb{F})$
given by $P\cdot A:=P^*AP$.
Note that for each matrix $A\in\Phi_{n, \mathbb{F}}^{-1}(0)$,
there exists $P\in U(n, \mathbb{F})$ such that
$$P\cdot A=\left\{
\begin{aligned}
\operatorname{diag}(\lambda_1, -\lambda_1, \cdots, \lambda_p, -\lambda_p, \lambda_{p+1})
\ \ \
&\mbox{if}\ n=2p+1\ \mbox{is odd};\\
\operatorname{diag}(\lambda_1, -\lambda_1, \cdots, \lambda_{p+1}, -\lambda_{p+1})
\ \ \ \ \,
&\mbox{if}\ n=2p+2\ \mbox{is even},
\end{aligned}\right.$$
where $\lambda_1\geq\cdots\geq\lambda_{p+1}\geq0$.
By analyzing the linear dependence of normal vectors
$\operatorname{grad}f_3, \cdots, \operatorname{grad}f_{2p+1}$,
we can show that
$A\in\mathcal{B}_{n, \mathbb{F}}$ if and only if $\lambda_1, \cdots, \lambda_{p+1}$ are distinct.
Thus we can characterize $\mathcal{B}_{n, \mathbb{F}}$ and $\mathcal{C}_{n, \mathbb{F}}$ more accurately.
It follows that
$$\begin{aligned}
&\mathcal{C}_{3, \mathbb{F}}=\varnothing,\\
&\mathcal{C}_{4, \mathbb{F}}=
U(n, \mathbb{F})\cdot\frac12\operatorname{diag}\(1, -1, 1, -1\),\\
&\mathcal{C}_{5, \mathbb{F}}=U(n, \mathbb{F})\cdot
\left\{\frac1{\sqrt2}\operatorname{diag}\(1, -1, 0, 0, 0\),
\frac12\operatorname{diag}\(1, -1, 1, -1, 0\)\right\}\\
&\mathcal{C}_{6, \mathbb{F}}=U(n, \mathbb{F})\cdot
\left\{\frac12\operatorname{diag}\(s, -s, s, -s, \sqrt{2(1-s^2)}, -\sqrt{2(1-s^2)}\)
: 0\leq s\leq 1\right\}.
\end{aligned}$$
For the case of $n=3$,
$\mathcal{B}_{3, \mathbb{R}}$, $\mathcal{B}_{3, \mathbb{C}}$ and $\mathcal{B}_{3, \mathbb{H}}$
are exactly Cartan's minimal isoparametric hypersurfaces with 3 distinct principal curvatures of multiplicity $m=1, 2$ and $4$, respectively.
For the case of $n\geq4$, $\mathcal{C}_{n, \mathbb{F}}\neq\varnothing$ implies that
$\mathcal{B}_{n, \mathbb{F}}$ is noncompact.
However, by finding principal directions at diagonal points, we can prove the following result:
\begin{theorem}{\rm(See \cite{GZ23})}\label{Bn austere}
$\mathcal{B}_{n, \mathbb{F}}$ is an orientable, connected, full, austere submanifold of $S(n, \mathbb{F})$ with flat normal bundle.
\end{theorem}

For the singular orbit $\mathcal{C}_{4, \mathbb{F}}$, we can show that
$\mathcal{C}_{4, \mathbb{F}}$ is a closed, full, austere submanifold
which is diffeomorphic to the Grassmann manifold $\operatorname{G}_2(\mathbb{F}^4)$.
Moreover, the Gauss equation yields that the scalar curvature of $\mathcal{B}_{4, \mathbb{R}}$ is
$$\rho(A)=42-2\(\operatorname{tr}A^4-\frac14\)^{-1}.$$
Since $\operatorname{tr}A^4=\frac14$ if and only if $A\in\mathcal{C}_{4, \mathbb{R}}$,
we know that $\mathcal{A}_4=\mathcal{B}_{4, \mathbb{R}}\cup \mathcal{C}_{4, \mathbb{R}}$ is not a smooth embedded submanifold of $S(4, \mathbb{R})$, unlike Cartan's isoparametric hypersurface $\mathcal{A}_3=\mathcal{B}_{3, \mathbb{R}}$ in the case of $n=3$.

Since for every austere subspace $\mathcal{Q}\subset\mathcal{S}_n$, we observe that
$\mathcal{Q}\cap S(n, \mathbb{R})$ is a totally geodesic sphere containing in
$\mathcal{A}_n=\mathcal{B}_{n, \mathbb{R}}\cup \mathcal{C}_{n, \mathbb{R}}$.
As an immediate application, we can obtain the following upper bound for the dimension of austere subspaces by the submanifold structure of $\mathcal{B}_{n, \mathbb{R}}$.
\begin{theorem}{\rm(See \cite{GZ23})}\label{dimest}
Let $\mathcal{Q}$ be an austere subspace in $M(n, \mathbb{R})$ such that
$\mathcal{Q}\cap\mathcal{B}_{n, \mathbb{R}}\neq\varnothing$, then
\begin{equation}\label{dimQ}
\operatorname{dim}\mathcal{Q}\leq\left\{\begin{aligned}
p^2+2p\ \ \ \ \ &\mbox{if}\ n=2p+1;\\
p^2+3p+2\ \ \ &\mbox{if}\ n=2p+2.
\end{aligned}\right.
\end{equation}
\end{theorem}

If $n=2p+2$, then the equality in (\ref{dimQ}) is achieved by the maximal austere subspace
$$\left\{\begin{pmatrix}
m_1 & m_2  \\
m_2 & -m_1 \\
\end{pmatrix}: m_1, m_2\in\mathcal{S}_{p+1}\right\}.$$
But for the case $n=2p+1$, we think the optimal upper bound should be smaller.
In addition, since $\mathcal{C}_{n, \mathbb{R}}$ is the union of singular orbits,
we believe that the intersection condition
$\mathcal{Q}\cap\mathcal{B}_{n, \mathbb{R}}\neq\varnothing$ in Theorem \ref{dimest} is not necessary.

\section*{3. An Application to Dupin Hypersurfaces}
We start this section by introducing some necessary concepts.
An oriented hypersurface $f: M^n \rightarrow N^{n+1}$ is called a Dupin hypersurface if along each curvature surface, the corresponding principal curvature is constant.
A Dupin hypersurface $M^n$ is called proper Dupin if
the number $g$ of distinct principal curvatures is constant.
In real space forms, it is easy to see that
every isoparametric hypersurface is a proper Dupin hypersurface.
In \cite{Pinkall}, Pinkall introduced four methods to construct a Dupin hypersurface from a lower dimensional Dupin hypersurface,
and he showed that there exists a proper Dupin hypersurface for any given number $g$ of distinct principal curvatures and any given multiplicities $m_1, \cdots, m_g$.
For compact proper Dupin hypersurfaces in spheres,
Thorbergsson \cite{Thorbergsson83} proved that the number $g$ of distinct principal curvatures can only be $1, 2, 3, 4$ or $6$, the same restriction as for isoparametric hypersurfaces proven by M\"{u}nzner \cite{Mun} (There is a new proof by Fang \cite{Fa17}).
See \cite{Cecil08, Cecil24, Chi19} for detailed surveys on the history and developments of isoparametric hypersurfaces and Dupin hypersurfaces.

In contrast to the situation for isoparametric hypersurfaces,
a proper Dupin hypersurface may not be extended to a compact proper Dupin hypersurface.
%Reducibility is a subtle but important property of Dupin submanifolds.
Hence reducibility is an important condition in the local classification of proper Dupin hypersurface.
A Dupin submanifold obtained from lower dimensional Dupin submanifold
via one of Pinkall's constructions is said to be reducible.
More generally, a Dupin submanifold which is locally Lie equivalent to
such a Dupin submanifold is said to be reducible \cite{Cecil}.
Cecil-Chi-Jensen \cite{CCJ} proved that every compact, connected proper Dupin
hypersurface in the Euclidean space with $g>2$ distinct principal curvatures is irreducible.
In the elegant survey \cite{Thorbergsson00}, Thorbergsson raised the following question:
\begin{ques0}
Is it possible that by assuming irreducibility instead of compactness of a proper Dupin
hypersurface the conclusion $g = 1, 2, 3, 4$ or $6$ can be drawn
as well as the restrictions on the multiplicities in Stolz's theorem \cite{Stolz}?
\end{ques0}

During the study of our new examples of austere submanifolds,
we find that there are three counterexamples for Thorbergsson's question.
In fact, we observed that there are $5$ distinct principal curvature functions
on the open submanifold
$$\widetilde{\mathcal{B}}_{4, \mathbb{F}}
=\mathcal{B}_{4, \mathbb{F}}\cap\operatorname{GL}(4, \mathbb{F})
=U(4, \mathbb{F})\cdot\widetilde{\mathcal{D}}_4,$$
where
$$\widetilde{\mathcal{D}}_4:=
\{\operatorname{diag}(\lambda_1, -\lambda_1, \lambda_2, -\lambda_2)\in\mathcal{B}_{4, \mathbb{F}}:
\lambda_1>\lambda_2>0\}.$$
Using some technologies from the theory of Lie sphere geometry \cite{Cecil},
we can prove the following result:
\begin{theorem}{\rm(See \cite{GZ23})}\label{B4_pD}
$\widetilde{\mathcal{B}}_{4, \mathbb{F}}$ is a $(6m+1)$-dimensional, orientable, connected, noncompact, irreducible proper Dupin hypersurface of the unit sphere $S(4, \mathbb{F})$ with $5$ distinct principal curvatures of multiplicities $m$, $m$, $m$, $m$ and $2m+1$, where $m=\operatorname{dim}_{\mathbb{R}}(\mathbb{F})$.
\end{theorem}

\section*{4. Normal Scalar Curvature Inequalities on a Class of Austere Submanifolds}

Austere submanifolds also arise from minimal Wintgen ideal submanifolds
which are minimal submanifolds attaining the following normal scalar curvature inequality (\ref{GeoDDVV}) everywhere \cite{DT09, XLMW14, XLMW18}.
In 1999, De Smet-Dillen-Verstraelen-Vrancken \cite{DDVV99} proposed the so-called DDVV conjecture:
For an immersed submanifold $M^n$ of a real space form $N^{n+m}(\kappa)$ with constant sectional curvature $\kappa$, there exists a pointwise normal scalar curvature inequality (also known as DDVV inequality)
\begin{equation}\label{GeoDDVV}
\rho+\rho^{\bot}\leq \|H\|^2+\kappa,
\end{equation}
where $\rho$ denotes the normalized scalar curvature, $\rho^{\bot}$ denotes the normalized normal scalar curvature and $\|H\|$ denotes the mean curvature.
%Dillen-Fastenakels-Veken \cite{DFV07} transformed (\ref{GeoDDVV}) into the so-called DDVV inequality
%\begin{equation}\label{AlgDDVV}
%\sum^m_{r,s=1}\|\[B_r,B_s\]\|^2\leq\(\sum^m_{r=1}\|B_r\|^2\)^2
%\end{equation}
%for $n\times n$ real symmetric matrices $B_1, \ldots, B_m$,
%where $[A,B]:=AB-BA$ denotes the commutator and
%$\|B\|^2:=\operatorname{tr}(BB^t)$ denotes the squared Frobenius norm.
The DDVV conjecture was proved by Lu \cite{Lu11} and Ge-Tang \cite{GT08} independently,
and the equality in (\ref{GeoDDVV}) holds at some point $p\in M$ if and only if there exist an orthonormal basis $\{e_1,\cdots,e_n\}$ of $T_pM$ and an orthonormal basis $\{\xi_1,\cdots,\xi_m\}$ of $T_p^{\bot}M$ such that
the matrices corresponding to the second fundamental form along directions $\xi_1,\cdots,\xi_m$ are
$$A_1=\lambda_1I_n+\mu\operatorname{diag}(H_1,0),\
A_2=\lambda_2I_n+\mu\operatorname{diag}(H_2,0),\
A_3=\lambda_3I_n$$ and $A_{r}=0$ for $r>3$,
where $\mu, \lambda_1, \lambda_2, \lambda_3$ are real numbers and $$H_1=
\begin{pmatrix}
1 & 0 \\
0 & -1 \\
\end{pmatrix},\
H_2=
\begin{pmatrix}
0 & 1 \\
1 & 0 \\
\end{pmatrix}.$$
See \cite{GLTZ24, GT11} for detailed surveys on the DDVV inequality.

For the focal submanifolds of isoparametric hypersurfaces in unit spheres,
Ge-Tang-Yan \cite{GTY20} proved sharper normal scalar curvature inequalities
and completely characterized the subsets achieving the equalities.
Recall that the isoparametric focal submanifolds are a rich source of austere submanifolds.
Thus it is natural to consider new normal scalar curvature inequalities on more general austere submanifolds.
For a given austere subspace $\mathcal{Q}$,
we say that an austere submanifold $M$ is of type $\mathcal{Q}$ if for each point $p\in M$,
there exists an orthonormal basis of $T_pM$ such that $|\operatorname{II}_p|\subset\mathcal{Q}$.
This definition, to some degree, generalizes the three types of $4$-dimensional austere submanifolds
introduced by Ionel-Ivey \cite{II10, II12}.
Since minimal Wintgen ideal submanifolds are austere submanifolds of type $\mathcal{Q}(n)$,
there is no new normal scalar curvature inequalities for this class of austere submanifolds.
However, for austere submanifolds of type $\mathcal{Q}(p, q)$,
we can establish the following sharper normal scalar curvature inequalities.

\begin{theorem}{\rm(See \cite{GTZ24})}\label{GeoDDVVQ(p,q)}
Let $M^n$ be an austere submanifold of type $\mathcal{Q}(p,q)$ in a real space form $N^{n+m}(\kappa)$. Then
\begin{equation}\label{GeoineqQ(p,q)}
\rho^{\bot}\leq \frac{\sqrt{2}}{2}(\kappa-\rho).
\end{equation}
The equality holds at some point $p\in M$ if and only if there exist an orthonormal basis $\{e_1,\cdots,e_n\}$ of $T_pM$ and an orthonormal basis $\{\xi_1,\cdots,\xi_m\}$ of $T_p^{\bot}M$ such that
the matrices corresponding to the second fundamental form along directions $\xi_1,\cdots,\xi_m$ are
$$A_1=\lambda\operatorname{diag}(C_1,0),\ A_2=\lambda\operatorname{diag}(C_2,0)\
\mbox{and}\ A_r=0\ \mbox{for}\ r>2,$$
where $\lambda$ is a real number and
$$C_1:=
\begin{pmatrix}
\ 0 & 0 & 1 & 0\ \\
\ 0 & 0 & 0 & 1\ \\
\ 1 & 0 & 0 & 0\ \\
\ 0 & 1 & 0 & 0\ \\
\end{pmatrix},\
C_2:=
\begin{pmatrix}
0 & 0 & 0 & -1\\
0 & 0 & 1  & 0\\
0 & 1 & 0 & 0\\
-1 & 0 & 0 & 0\\
\end{pmatrix}.$$
\end{theorem}

The proof follows Ge-Tang's method \cite{GT08} which enables us to derive the equality condition.
Except for the trivial case $A_1=\cdots=A_m=0$,
the equality in (\ref{GeoineqQ(p,q)}) can only be achieved when $p, q\geq 2$.
For the remaining cases, we have the following sharper inequality.

\begin{theorem}{\rm(See \cite{GTZ24})}\label{GeoDDVVQ(n-1,1)}
Let $M^n$ be an austere submanifold of type $\mathcal{Q}(n-1,1)$ in a real space form $N^{n+m}(\kappa)$. Then for any point $p\in M$,
\begin{equation}\label{GeoineqQ(n-1,1)}
\rho^{\bot}\leq \sqrt{\frac{d-1}{2d}}(\kappa-\rho),
\end{equation}
where $d$ denotes the dimension of $|\operatorname{II}_p|$.
The equality holds if and only if there exist an orthonormal basis $\{e_1,\cdots,e_n\}$ of $T_pM$ and an orthonormal basis $\{\xi_1,\cdots,\xi_m\}$ of $T_p^{\bot}M$ such that
the matrices corresponding to the second fundamental form along directions $\xi_1,\cdots,\xi_m$ are
$$A_r=\lambda \(E_{rn}+E_{nr}\)\ \mbox{for}\ 1\leq r\leq d\
\mbox{and}\ A_r=0\ \mbox{for}\ r>d,$$
where $\lambda$ is a real number and
$E_{ij}$ denotes the $n\times n$ matrix with $1$ in position $(i,j)$ and $0$ elsewhere.
\end{theorem}

To illustrate the above theorems make sense,
we give some explicit examples of austere submanifold of type $\mathcal{Q}(p,q)$ in \cite{GTZ24}.
It is worth pointing out that the isoparametric focal submanifold $M_+$
with $g=4$ and multiplicities $(1, 2)$ achieves the equality in (\ref{GeoineqQ(p,q)}) everywhere, and
some of Bryant's generalized helicoids achieve the equality in (\ref{GeoineqQ(n-1,1)}) on certain subsets.
As a byproduct, we can establish the following Simons-type gap theorem for closed austere submanifolds of type $\mathcal{Q}(p,q)$ in unit spheres.

\begin{theorem}{\rm(See \cite{GTZ24})}\label{Simons}
Suppose $M^n$ is a closed, connected, austere submanifold of type $\mathcal{Q}(p,q)$ in the unit sphere $\mathbf{S}^{n+m}$.
Then $$\int_M S(S-n)\,dV_M\geq0.$$
In addition, if $S\leq n$, then $M$ is either a totally geodesic sphere or
the Clifford torus $S^p\(\sqrt{\frac12}\)\times S^p\(\sqrt{\frac12}\)$ with $n=2p$.
\end{theorem}

It is natural to compare Theorem \ref{Simons} with the famous rigidity theorem:
\begin{theorem}{\rm(See \cite{CDK, Lawson, Simons})}
Suppose $M^n$ is a closed, connected, minimal submanifold in the unit sphere $\mathbf{S}^{n+m}$.
Then $$\int_M S\[\(2-\frac{1}{m}\)S-n\]\,dV_M\geq0.$$
In addition, if $S\leq n/\(2-\frac{1}{m}\)$, then $M$ is either a totally geodesic sphere or
a Clifford torus $S^k\Big(\sqrt{\frac{k}{n}}\Big)\times S^{n-k}\Big(\sqrt{\frac{n-k}{n}}\Big)$
or the Veronese surface $\mathbb{R}P^2\subset\mathbf{S}^{4}$,
where $1\leq k\leq n-1$.
\end{theorem}

We can see that when replacing minimal submanifolds with austere submanifolds of type $\mathcal{Q}(p,q)$
in the assumption, the first gap of $S$ is significantly larger, and the similar rigidity result still holds.

%\newpage
\markboth{\small {\sc Jianquan Ge \& Yi Zhou}}{\small Recent Progress on Austere Submanifolds}


\begin{thebibliography}{HHHH}

\markboth{\small {\sc Jianquan Ge \& Yi Zhou}}{\small Recent Progress on Austere Submanifolds}

%\bibitem{AKMR09}
%K. M. R. Audenaert, \emph{Variance bounds, with an application to norm bounds for commutators},
%Linear Algebra Appl. \textbf{432} (5) (2009), 1126-1143.

%\bibitem{BCO}
%J. Berndt, S. Console and C. E. Olmos, \emph{Submanifolds and holonomy. Second edition}. Monographs and Research Notes in Mathematics. CRC Press, Boca Raton, FL, 2016.

\bibitem{Bryant}
R. Bryant, \emph{Some Remarks on the Geometry of Austere Manifolds}, Bol. Soc. Brasil. Mat. (N.S.)
\textbf{21} (1991), 122-157.

%\bibitem{BW05}
%A. B\"{o}ttcher and D. Wenzel, \emph{How big can the commutator of two matrices be and how big is it typically?} Linear Algebra Appl. \textbf{403} (2005), 216-228.

%\bibitem{BW08}
%A. B\"{o}ttcher and D. Wenzel, \emph{The Frobenius norm and the commutator}, Linear Algebra Appl. \textbf{429} (8-9) (2008), 1864-1885.

\bibitem{Cartan1}
E. Cartan, \emph{Sur des familles remarquables d'hypersurfaces isoparam\'{e}tri-
ques dans les espaces sph\'{e}riques}, Math. Z. \textbf{45} (1939), 335-367.

\bibitem{Cartan2}
E. Cartan, \emph{Sur quelques familles remarquables d'hypersurfaces},
C. R. Congr\`{e}s Math. Li\`{e}ge (1939), 30-41.

\bibitem{Cecil08}
T. E. Cecil, \emph{Isoparametric and Dupin hypersurfaces},
SIGMA Symmetry Integrability Geom. Methods Appl. \textbf{4} (2008), Paper 062.

\bibitem{Cecil}
T. E. Cecil, \emph{Lie Sphere Geometry}, 2nd edn. Springer, New York, 2008.

\bibitem{Cecil24}
T. E. Cecil, \emph{Classifications of Dupin Hypersurfaces in Lie Sphere Geometry}, Acta Math. Sci. Ser. B (Engl. Ed.) \textbf{44} (2024), no. 1, 1-36.

\bibitem{CCJ}
T.E. Cecil, Q. S. Chi and G.R. Jensen, \emph{Dupin hypersurfaces with four principal curvatures II}, Geom. Dedicata \textbf{128} (2007), 55-95.

%\bibitem{CCJ07}
%T. E. Cecil, Q. S. Chi, and G. R. Jensen, \emph{Isoparametric hypersurfaces with four principal curvatures}, Ann. of Math. (2) \textbf{166} (2007), no. 1, 1-76.

\bibitem{CR85}
T. E. Cecil and P. J. Ryan, \emph{Tight and Taut Immersions of Manifolds}, Research Notes in Mathematics, vol. 107, Pitman, London, 1985.

\bibitem{CR15}
T. E. Cecil and P. J. Ryan, \emph{Geometry of hypersurfaces}, Springer Monographs in Mathematics,
Springer, New York, 2015.

%\bibitem{CVW10}
%C. Cheng, S. Vong and D. Wenzel, \emph{Commutators with maximal Frobenius norm}, Linear Algebra Appl.
%\textbf{432} (2010), 292-306.

\bibitem{CDK}
S. S. Chern, M. do Carmo and S. Kobayashi,
\emph{Minimal submanifolds of a sphere with second fundamental form of constant length},
in: Functional Analysis and Related Fields, Proc. Conf. for M. Stone, Univ. Chicago, Chicago, Ill., 1968, Springer, New York, 1970, pp. 59-75.

%\bibitem{Chi11}
%Q. S. Chi, \emph{Isoparametric hypersurfaces with four principal curvatures, II}, Nagoya Math. J. \textbf{204} (2011), 1-18.

%\bibitem{Chi13}
%Q. S. Chi, \emph{Isoparametric hypersurfaces with four principal curvatures, III}, J. Differential Geom. \textbf{94} (2013), 487-522.

%\bibitem{Chi16}
%Q. S. Chi, \emph{Isoparametric hypersurfaces with four principal curvatures, IV},  J. Differential Geom. \textbf{115}(2) (2020),  225-301.

\bibitem{Chi19}
Q. S. Chi, \emph{The isoparametric story, a heritage of \'{E}lie Cartan}, Proceedings of the International Consortium of Chinese Mathematicians, 2018, International Press of Boston (2020), 197-260.

\bibitem{CL08}
T, Choi and Z. Lu, \emph{On the DDVV conjecture and the comass in calibrated geometry (I)}, Math. Z. \textbf{260} (2008), 409-429.

\bibitem{DF01}
M. Dajczer and L. A. Florit, \emph{A class of austere submanifolds}, Illinois J. Math. \textbf{45} (2001), no. 3, 735-755.

\bibitem{DFT}
M. Dajczer, L. Florit and R. Tojeiro, \emph{Reducibility of Dupin submanifolds}, Illinois J. Math. \textbf{49} (2005), 759-791.

\bibitem{DT09}
M. Dajczer and R. Tojeiro, \emph{Submanifolds of codimension two attaining equality in an extrinsic inequality}, Math. Proc. Cambridge Philos. Soc. \textbf{146} (2009), no. 2, 461-474.

\bibitem{DV23}
M. Dajczer and T. Vlachos, \emph{A new class of austere submanifolds}, Comm. Anal. Geom. \textbf{31} (2023), no.4, 879-893.

\bibitem{DDVV99}
P. J. De Smet, F. Dillen, L. Verstraelen and L. Vrancken, \emph{A pointwise inequality in submanifold
theory}, Arch. Math. (Brno), \textbf{35} (1999), 115-128.

%\bibitem{DFV07}
%F. Dillen, J. Fastenakels and J. Veken,
%\emph{Remarks on an inequality involving the normal scalar curvature},
%Proceedings of the International Congress on Pure and Applied Differential Geometry
%PADGE, Brussels, edited by F. Dillen and I. Van deWoestyne, Shaker, Aachen, 2007: 83-92.

%\bibitem{DHTV04}
%F. Dillen, S. Haesen, M. Petrovi\'{c}-Torga\v{s}ev and L. Verstraelen, \emph{An inequality between
%intrinsic and extrinsic scalar curvature invariants for codimension 2 embeddings}, J. Geom. Physics
%\textbf{52} (2004), 101-112.

%\bibitem{DN85}
%J. Dorfmeister and E. Neher, \emph{Isoparametric hypersurfaces, case $g = 6, m =1$}, Comm. Algebra \textbf{13} (1985), 2299-2368.

\bibitem{Fa17}
F. Q. Fang, \emph{Dual submanifolds in rational homology spheres}, Sci. China Math. \textbf{60} (2017), no. 9, 1549-1560.

%\bibitem{FKM81}
%D. Ferus, H. Karcher and H. F. M\"{u}nzner, \emph{Cliffordalgebren und neue isoparametrische Hyperfl\"{a}chen}, Math. Z. \textbf{177} (1981), 479-502.

%\bibitem{Ge14}
%J. Q. Ge, \emph{DDVV-type inequality for skew-symmetric matrices and Simons-type inequality for Riemannian submersions}, Adv. Math. \textbf{251} (2014), 62-86.

%\bibitem{GLZ21}
%J. Q. Ge, F. G. Li and Y. Zhou, \emph{Some generalizations of the DDVV and BW inequalities},
%Trans. Amer. Math. Soc. \textbf{374} (2021), no. 8, 5331-5348.

%\bibitem{GLLZ20}
%J. Q. Ge, F. G. Li, Z. Lu and Y. Zhou, \emph{On some conjectures by Lu and Wenzel},
%Linear Algebra Appl. \textbf{592} (2020), 134-164.

\bibitem{GLTZ24}
J. Q. Ge, F. G. Li, Z. Z. Tang and Y. Zhou, \emph{A Survey on the DDVV-type Inequalities}, Adv. Math. (China) \textbf{53} (2024), no. 3, 449-467.

\bibitem{GT08}
J. Q. Ge and Z. Z. Tang, \emph{A proof of the DDVV conjecture and its equality case}, Pacific J. Math. \textbf{237} (2008), no. 1, 87-95.

\bibitem{GT11}
J. Q. Ge and Z. Z. Tang, \emph{A survey on the DDVV conjecture}, \textit{Harmonic maps and differential geometry}, 247-254, Contemp. Math. \textbf{542}, Amer. Math. Soc., Providence, RI, 2011.

\bibitem{GT22}
J. Q. Ge and Z. Z. Tang, \emph{Isoparametric Polynomials and Sums of Squares}, Int. Math. Res. Not. IMRN 2023, no. 24, 21226-21271.

\bibitem{GTY20}
J. Q. Ge, Z. Z. Tang and W. J. Yan, \emph{Normal scalar curvature inequality on the focal submanifolds of isoparametric hypersurfaces}, Int. Math. Res. Not. IMRN 2020, no. 2, 422-465.

\bibitem{GTZ24}
J. Q. Ge, Y. Tao and Y. Zhou, \emph{Normal scalar curvature inequality on a class of austere submanifolds}, arXiv:2410.17761.

%\bibitem{GXYZ17}
%J. Q. Ge, S. Xu, H. Y. You and Y. Zhou, \emph{DDVV-type inequality for Hermitian matrices}, Linear Algebra Appl. \textbf{529} (2017), 133-147.

\bibitem{GZ23}
J. Q. Ge and Y. Zhou, \emph{Austere matrices, austere submanifolds and Dupin hypersurfaces}, arXiv:2302.06105.

\bibitem{HL}
R. Harvey and H. B. Lawson, \emph{Calibrated Geometries}, Acta Math. \textbf{148} (1982), 47-157.

\bibitem{II10}
M. Ionel and T. Ivey, \emph{Austere submanifolds of dimension four: examples and maximal types}, Illinois J. Math. \textbf{54} (2010), no. 2, 713-746.

\bibitem{II12}
M. Ionel and T. Ivey, \emph{Ruled austere submanifolds of dimension four}, Differential Geom. Appl. \textbf{30} (2012), no. 6, 588-603.
%\bibitem{II16}
%M. Ionel and T. A. Ivey, \emph{Austere submanifolds in \mathbb{C}\mathbf{P}^n}, Comm. Anal. Geom. \textbf{24} (2016), no. 4, 821-841.

\bibitem{IST09}
O. Ikawa, T. Sakai and H. Tasaki, \emph{Weakly reflective submanifolds and austere submanifolds},
J. Math. Soc. Japan \textbf{61} (2009), no. 2, 437-481.

%\bibitem{JM84}
%L. P. Jorge and F. Mercuri, \emph{Minimal immersions into space forms with two principal curvatures},
%Math. Z. \textbf{187} (1984), no.3, 325-333.

\bibitem{KM22}
T. Kimura and K. Mashimo, \emph{Classification of Cartan embeddings which are austere submanifolds}, Hokkaido Math. J. \textbf{51} (2022), no. 1, 1-23.

\bibitem{Lawson}
H. B. Lawson, \emph{Local rigidity theorems for minimal hypersurfaces}, Ann. of Math. (2) \textbf{89} (1969), 187-197.

%\bibitem{Lu07}
%Z. Lu, \emph{Recent developments of the DDVV conjecture}, Bull. Transilv. Univ. Bra\c{s}ov Ser. B (N.S.) \textbf{14}(49) (2007), 133-143.

\bibitem{Lu11}
Z. Lu, \emph{Normal scalar curvature conjecture and its applications}, J. Funct. Anal. \textbf{261} (2011), 1284-1308.

%\bibitem{Lu12}
%Z. Lu, \emph{Remarks on the B\"{o}ttcher-Wenzel inequality}, Linear Algebra Appl. \textbf{436} (2012), no. 7, 2531-2535.

%\bibitem{MH73}
%J. Milnor and D. Husemoller, \emph{Symmetric Bilinear Forms}, Ergebnisse der Mathematik und ihrer Grenzgebiete, Band 73, Springer-Verlag, New York-Heidelberg, 1973.

%\bibitem{Miy13}
%R. Miyaoka, \emph{Isoparametric hypersurfaces with $(g,m)=(6,2)$}, Ann. of Math. (2) \textbf{177} (2013), 53-110.

%\bibitem{Miy16}
%R. Miyaoka, \emph{Errata of ``Isoparametric hypersurfaces with $(g,m)=(6,2)$"}, Ann. of Math. (2) \textbf{183} (2016), 1057-1071.

\bibitem{Mun}
H. F. M{\"u}nzner, \emph{Isoparametrische Hyperfl{\"a}chen in Sph{\"a}ren, I and II}, Math. Ann. \textbf{251} (1980), 57--71 and \textbf{256} (1981), 215-232.

%\bibitem{OT75}
%H. Ozeki and M. Takeuchi. \emph{On some types of isoparametric hypersurfaces in spheres I},
%Tohoku Math. J. \textbf{27} (1975), 515-59.

\bibitem{Pinkall}
U. Pinkall, \emph{Dupin hypersurfaces}, Math. Ann. \textbf{270} (1985), 427-440.

%\bibitem{Scharlau}
%W. Scharlau, \emph{Quadratic and Hermitian forms}, Grundlehren der mathematischen Wissenschaften, Springer-Verlag, Berlin, 1985.

\bibitem{Simons}
J. Simons, \emph{Minimal varieties in Riemannian manifolds}, Ann. of Math. (2) \textbf{88} (1968), 62-105.

\bibitem{Stolz}
S. Stolz, \emph{Multiplicities of Dupin hypersurfaces}, Invent. Math. \textbf{138} (1999), no. 2, 253-279.

%\bibitem{TZ20}
%Z. Z. Tang and Y. S. Zhang, \emph{Minimizing cones associated with isoparametric foliations}, J. Differential Geom. \textbf{115} (2020), no. 2, 367-393.

\bibitem{Thorbergsson83}
G. Thorbergsson, \emph{Dupin hypersurfaces}, Bull. London Math. Soc. \textbf{15} (1983), 493-498.

\bibitem{Thorbergsson00}
G. Thorbergsson, \emph{A survey on isoparametric hypersurfaces and their generalizations}, Handbook of differential geometry, Vol. I, 963-995, North-Holland, Amsterdam, 2000.

%\bibitem{VJ08}
%S. W. Vong and X. Q. Jin, \emph{Proof of B\"{o}ttcher and Wenzel's conjecture}, Oper. Matrices \textbf{2} (2008), no. 3, 435-442.

\bibitem{XLMW14}
Z. X. Xie, T. Z. Li, X. Ma and C. P. Wang, \emph{M\"{o}bius geometry of three-dimensional Wintgen ideal submanifolds in $\mathbb{S}^5$}, Sci. China Math. \textbf{57} (2014), no. 6, 1203-1220.

\bibitem{XLMW18}
Z. X. Xie, T. Z. Li, X. Ma and C. P. Wang, \emph{Wintgen ideal submanifolds: reduction theorems and a coarse classification}, Ann. Global Anal. Geom. \textbf{53} (3) (2018), 377-403.

\end{thebibliography}
\end{document}